\documentclass[10pt]{amsart}
\usepackage{amsmath,amsthm,amssymb,color}
\usepackage{epsfig,amsfonts,latexsym, hyperref, bbm, mathrsfs}
\usepackage{natbib}
\usepackage{geometry}
\numberwithin{equation}{section}

\newcommand{\scrm}{\mathscr{M}}

\newcommand{\inv}{{-1}}
\newcommand{\bbs}{\mathbb{S}}

\newcommand{\real}{\mathbb{R}}

\newcommand{\stas}{\mbox{St$\alpha$S}}

\newcommand{\beq}{\begin{equation}}
\newcommand{\eeq}{\end{equation}}
\newcommand{\alns}[1]{\begin{align*}#1\end{align*}}
\newcommand{\aln}[1]{\begin{align} #1 \end{align}}
\newcommand{\been}{\begin{enumerate}}
\newcommand{\een}{\end{enumerate}}

\newcommand{\bdelta}{\boldsymbol{\delta}}

\newcommand{\eqd}{\,{\buildrel d \over =}\,}
\newcommand{\ckr}{C_c^+ (\bar{\real}_0)}
\newcommand{\ckrfty}{C_c^+(\rmfty)}
\newcommand{\zrer}{\bar{\real}_0}

\newcommand{\point}{\mathcal{P}}
\newcommand{\malpha}{m_\alpha}

\newcommand{\bbo}{\mathbbm{1}}
\newcommand{\cmax}{c_{max}}
\newcommand{\Novery}{N^{(y)}}

\newcommand{\Nundery}{N_{(y)}}

\newcommand{\rmfty}{\bar{\real}_{-\infty}}
\newcommand{\nmfty}{\mathscr{M}_{-\infty}}

\newcommand{\npreal}{\mathscr{M}_+}

\newcommand{\scrt}{\mathscr{T}}
\newcommand{\ckrp}{C_c^+(\real_+)}

\DeclareMathOperator{\maxmod}{maxmod}
\DeclareMathOperator{\fralpha}{\Phi_\alpha}

\DeclareMathOperator{\suplevel}{SL}
\DeclareMathOperator{\Exp}{\mathbf{Exp}}
\DeclareMathOperator{\Log}{\mathbf{Log}}
\DeclareMathOperator{\dtv}{d}

\usepackage{cleveref}
\newtheorem{thm}{Theorem}[section]

\newtheorem{cor}[thm]{Corollary}

\theoremstyle{remark}

\theoremstyle{definition}
\newtheorem{defn}[thm]{Definition}
\newtheorem{fact}[thm]{Fact}

\DeclareMathOperator{\prob}{\mathbf{P}}
\DeclareMathOperator{\exptn}{\mathbf{E}}

\DeclareMathOperator{\mbfs}{\mathbf{S}}

\newcommand{\eql}{\stackrel{d}{=}}

\allowdisplaybreaks

\hypersetup{
    colorlinks=true,
    linkcolor=red,
    citecolor=blue,
    urlcolor=blue,
    pdfborder={0 0 0}
}

\title[randomly scaled scale-decorated Poisson point process]{A note on randomly scaled scale-decorated Poisson point processes}

\author{ Ayan Bhattacharya }

\date{\today }

\begin{document}

\maketitle

\begin{abstract}
 Randomly scaled scale-decorated Poisson point process is introduced recently in \cite{bhattacharya:hazra:roy:2014} where it appeared as weak limit of a sequence of point processes in the context of branching random walk. In this article, we obtain a characterization for these processes based on scaled-Laplace functional. As a consequence, we obtain a characterization for strictly $\alpha$-stable point process (also known as scale-decorated Poisson point process) based on scaled-Laplace functional . a connection with randomly shifted decorated Poisson point process is obtained. The tools and approach used e very similar to those in \cite{subag:zeitouni:2015}.
\end{abstract}



\section{Introduction}

In the context of extreme value theory for independently and identically distributed random variables, Poisson point process appear as the limiting extremal process (see Chapter~3 in \cite{resnick:1987}, Chapter~7 in  \cite{resnick:2007} and references therein). Recently, it has been observed that the dependence among the random variables causes slower rate of growth for extremes (e.g. \cite{samorodnitsky:2004a}, \cite{roy:samorodnitsky:2008}, \cite{roy:2010}), clusters around extreme points (\cite{davis:hsing:1995}).  So the class of all Poisson point processes is not large enough to accommodate all possible extremal processes which motivates and necessitates different generalization of Poisson point processes. One such generalization is randomly scaled scale-decorated Poisson point process (SScDPPP) introduced in \cite{bhattacharya:hazra:roy:2014} in the context of branching random walk with displacements having regularly varying tail. In this article, we shall investigate properties of SScDPPP and derive a characterization based on its scaled-Laplace functional.

Consider a Poisson point process $\Lambda = (\lambda_i : \i \ge 1)$ on $(0, \infty)$ and a collection of independently and identically distributed point processes $(\mathcal{P}_i : i \ge 1)$ on $\mathbb{R}$ which is independent of $\Lambda$. Informally speaking, we multiply the points of $\mathcal{P}_i$ by the $i$th Poisson point $\lambda_i$ for every $i \ge 1$ and superpose them. Then the superposition (may not always exist) is called scale-decorated Poisson point process (ScDPPP). If we multiply each point of ScDPPP by a positive random variable $U$ which is independent of $\Lambda$ and collection $(\mathcal{P}_i : i \ge 1)$, then we obtain SScDPPP. The notion ScDPPP is introduced in \cite{davydov:molchanov:zuyev:2008} as LePage series representation of a strictly $\alpha$-stable ($\stas$) point process. In this article, we shall introduce the notion scaled-Laplace functional to obtain characterization of SScDPPP. As a consequence, we obtain a characterization for ScDPPP based on scaled-Laplace functional. This characterization also has been used in \cite{bhattacharya:hazra:roy:2016} and \cite{bhattacharya:maulik:palmowski:roy:2016} to verify that the limiting extremal process is SScDPPP in the context of branching random walk with displacements having regularly varying tail. In the former article, this characterization has also been used to establish that the approximately scaled superposition of regularly varying point processes converges to $\stas$ point process. The tools and approach used in this article are motivated from \cite{subag:zeitouni:2014} and very similar to that article.

Note that the analogue of ScDPPP and SScDPPP in the context of random variables with exponentially decaying tail are decorated Poisson point process (DPPP) and randomly shifted decorated Poisson point process (SDPPP). In \cite{brunet:derrida:2011}, the notion DPPP is introduced in the context of branching random walk with exponentially decaying tail and branching Brownian motion. Then  characterization for DPPP is obtained in \cite{maillard:2013} based on shift-Laplace functional  which has been used in \cite{madaule:2011} to derive weak limit of appropriately shifted point processes. In \cite{subag:zeitouni:2014}, the notion SDPPP is introduced and a characterization is obtained based on shift-Laplace functional.

In this article, we also derived a connection between SScDPP and SDPPP. We shall show that every SDPPP on $(-\infty, \infty)$ can be transformed to a SCDPPP on $(0, \infty)$ and every ScDPPP on $(0, \infty)$ can be transformed to SDPPP on $(-\infty, \infty)$.

The article is organized as follows. In Section~\ref{sec:notation}, we have introduced some notations and stated main results of this article. The expression for scaled Laplace functional of SScDPPP has been derived in Section~\ref{sec:slf}.  In Section~\ref{sec:prop_scaledLaplace}, the properties of the scaled Laplace functional has been studied. The construction of the scale-decoration is given in Section~\ref{sec:scale_decoration}.  In Section~\ref{sec:rest_proofs},  proofs of main results and the connection between SSDPPP and SScDPPP are given.

\section{Notation and main results} \label{sec:notation}

\subsection{Notation}


Suppose $\zrer$ denote the punctured space $[-\infty, \infty] \setminus \{0\}$ and $\scrm(\zrer)$ denote the space of all Radon measures on $\zrer$ which does not put any mass on $\{\pm \infty\}$. The scalar multiplication on $\scrm(\zrer)$ by a positive real number $b$ is denoted by $\mbfs_b$ and is defined as follows: if $\point = \sum_i \bdelta_{u_i} \in \scrm(\zrer)$, then
\alns{
\mbfs_b \point = \sum_i \bdelta_{b u_i}. 
} 

In other words, a scalar multiple of a point measure is obtained by multiplying each point of the measure by a positive real number. A point process on $\zrer$ is an $\scrm(\zrer)$-valued random variable defined on $(\Omega, \mathcal{F}, \prob)$ that does not charge any mass at $\pm \infty$. The following definition of strictly $\alpha$-stable point process was introduced in \cite{davydov:molchanov:zuyev:2008}.

\begin{defn}[{\bf $\stas$ point process}; \cite{davydov:molchanov:zuyev:2008}] A point process $N$ (on $\zrer$) is called a strictly $\alpha$-stable ($\stas$) \label{not:stas} point process ($\alpha >0$) if for every $b_1, b_2 >0$,
\aln{
\mbfs_{b_1} N_1 + \mbfs_{b_2} N_2 \eqd \mbfs_{(b_1^\alpha + b_2^\alpha)^{1/\alpha}} N, \label{eq:chap2_alphastablepp}
}
where $N_1, N_2$ are independent copies of $N$, $+$ denotes superposition of point processes and $\eqd$ denotes equality in distribution.
\end{defn}

A point process $M$ will be called \textbf{randomly scaled strictly $\alpha$-stable point process} if $M \eqd \mbfs_U N$ where $N$ is a $\stas$ point process and $U$ is a positive random variable independent of $N$.

It has been established in \cite{davydov:molchanov:zuyev:2008}, that a point process (on $\real$) is $\stas$ if and only if it admits a series representation of a special kind (analogous to the LePage series representation for stable processes). To be more precise, we need the following definition which is introduced in \cite{bhattacharya:hazra:roy:2014}.

\begin{defn}[Scale-decorated Poisson point process]
A point process $N$ is called a scale decorated Poisson point process with intensity measure $m$ and scale-decoration $\point$ (denoted by $N \sim ScDPPP(m, \point)$) if there exists a Poisson random measure $\Lambda = \sum_{i=1}^\infty \bdelta_{\lambda_i}$ on $(0, \infty)$ with intensity measure $m$ and a point process $\point$ such that 
\alns{
N \eqd \sum_{i=1}^\infty \mbfs_{\lambda_i} \point
}
where $(\point_i : i \ge 1)$ is a collection of independent copies of $\point$.
\end{defn}

As mentioned above, it has been observed in \cite{davydov:molchanov:zuyev:2008} (see Example~8.6 therein) that a point process $N$ is $\stas$ if and only if $N \sim ScDPPP(m_\alpha, \point)$ where $\point$ is a point process on $\zrer$ and $m_\alpha$ is  measure on $(0, \infty)$ with $m_\alpha\Big((x, \infty) \Big) = \alpha x^{-\alpha -1}$ for every $x>0$. The light-tailed analogue of this result has been proved in a novel approach by \cite{maillard:2013}. 

A point process $M$ is called a \textbf{randomly scaled scale-decorated Poisson point process} with mean measure $m$, scale-decoration $\point$ and random scale $U$ (denoted by $M \sim SScDPPP(m, \point, U)$) if $M \eqd \mbfs_U N$ where $N \sim ScDPPP(m, \point)$ and $U$ is a positive random variable independent of $N$.

 Let $\ckr$ \label{not:ckr} denote the space of all nonnegative continuous functions defined on $\bar{\mathbb{R}}_0$ with compact support (and hence vanishing in a neighbourhood of $0$). By an abuse of notation, for a measurable function $f: \bar{\mathbb{R}}_0 \to [0, \infty)$, we denote by $\mbfs_{y}f(\cdot)$ the function $f(y \cdot)$. For a point process $N$ on $\bar{\real}_0$ and any $y>0$,  one has $\int f \dtv \mbfs_y N = \int \mbfs_y f \dtv N$. The Laplace functional of a point process $N$ will be denoted by
\aln{
\Psi_N  (f) = \exptn \Big( \exp \Big\{- N(f) \Big\} \Big), \label{eq:chap2_lap_not}
}
where $N(f) = \int f \dtv N$. \label{not:int_measure} In parallel to the notion of shift-Laplace functional from \cite{subag:zeitouni:2014}, we define the {\bf scaled Laplace functional} as 
\aln{
\Psi_N(f\|y) := \Psi_N(\mbfs_{y^\inv} f)  \label{not:chap2_scale_laplace} 
}
for some $y>0$. Let $g: (0,\infty) \to (0, \infty)$ be a measurable function. We define $[g]_{sc}$ as the class of all positive measurable functions $f: \mathbb{R}_+ \to \mathbb{R}_+$ such that for some $y>0$, $f(x) = g(yx)$ for all $x >0$. Let us define by $\fralpha(x)$ the Frech\'et distribution function, i.e., for each $\alpha>0$,
\aln{
 \fralpha(x)=\exp(-x^{-\alpha}), \qquad x>0. \label{not:chap2_fralpha}
 }

\begin{defn}[Scale-uniquely supported]
The scaled Laplace functional of the point process $N$ is uniquely supported on $[g]_{sc}$ if  for any $f \in \ckr$, there exists a constant $c_f$ (depending on $f$ only) such that $\Psi_N( f\|y )= g(y c_f)$ for all $y>0$.
\end{defn}

\subsection{Main Result}

The notion of scale-uniquely supported is intimately tied to the behaviour of SScDPPP. Sometimes, it is not possible to write down the SScDPPP representation, but it is easy to study properties of its scaled Laplace functional. In those cases, the following theorem turns out to be useful; see e.g., \cite{bhattacharya:hazra:roy:2014} and \cite{bhattacharya:maulik:palmowski:roy:2016}. The theorem can be very useful to study the weak limit of a sequence of  point processes as existence of the weak limit can be guaranteed by studying properties of scaled Laplace functional of the limit.

\begin{thm}\label{thm:sub_zet}
Let $N$ be a locally finite point process on $\zrer$ satisfying the following assumptions:
\aln{
\prob(N(\bar{\real}_0) > 0) >0 \mbox{ and  }  \exptn \Big( N(\bar{\real}_0 \setminus (-a,a)) \Big) < \infty \label{eq:chap2_ass_sub_zet}
 }
for some $a>0$. Then the following statements are equivalent:
\let\myenumi\theenumi
\let\mylabelenumi\labelenumi
\renewcommand{\theenumi}{Prop\myenumi}
\renewcommand{\labelenumi}{{\rm (\theenumi)}}
\begin{enumerate}

\item $\Psi_N(f\|\cdot)$ is scale-uniquely supported on $[g]_{sc}$ for all $f\in \ckr$ for some function $g : \real_+ \to \real_+$. \label{sub:zet:prop1}

\item $\Psi_N(f\|\cdot)$ is scale-uniquely supported on $[g]_{sc}$  for all $f\in \ckr$, where
\beq \label{eq:chap2_SLF}
g(y) = \exptn \Big( \fralpha( ycW^\inv)  \Big),
\eeq
for some $\alpha >0$, for some $c >0$ and some positive random variable $W$. \label{sub:zet:prop2}

\item $N \sim SScDPPP( \malpha, \point,W)$ for some point process $\point$, some positive random variable $W$ and some positive scalar $\alpha >0$ (same as in \ref{sub:zet:prop2})  where $m_\alpha(\cdot)$ is a measure on $(0,\infty)$, such that $m_\alpha ((x, \infty)) = x^{-\alpha}$ for every $x >0$. \label{sub:zet:prop3}

\end{enumerate}

\end{thm}

 The next result is an immediate corollary of the above proposition.  In \cite{davydov:molchanov:zuyev:2008}, it has been established that  \ref{item:DM1} and \ref{item:DM2} are equivalent.
\begin{cor} \label{propn:strict_stable}
Assume that $\prob(N(\zrer) >0) >0$. Fix $\alpha >0$. Then the following statements are equivalent:
\begin{enumerate}
\let\myenumi\theenumi
\let\mylabelenumi\labelenumi
\renewcommand{\theenumi}{B\myenumi}
\renewcommand{\labelenumi}{{\rm (\theenumi)}}
\item $N$ is a scale-decorated Poisson point process on $\zrer$ with Poisson intensity $m_\alpha(dx)$, where $m_\alpha$ is as defined in~\ref{sub:zet:prop3}. \label{item:DM1}

\item $N$ is a strictly $\alpha$-stable point process.\label{item:DM2}

\item The scaled Laplace functional of $N$ is scale-uniquely supported on the class $[\fralpha]_{sc}$.\label{item:DM3}

\end{enumerate}
\end{cor}

Suppose that  the  point process $N$ in Theorem~\ref{thm:sub_zet}  is supported on the positive part of the real line. Then these results are equivalent to the results in \cite{subag:zeitouni:2014} via  one to one correspondence between the spaces $\mathscr{M}((0, \infty])$ and $\mathscr{M}((-\infty,\infty])$ given by $\sum {\delta_{a_i}} \leftrightarrow \sum {\delta_{\log{a_i}}}$. In particular, the assumption of monotonicity of $g$ can be dropped from Corollary~3 of the aforementioned reference.

\subsection{Connection to SDPPP}

Here we shall establish a connection between the notions SScDPPP and SDPPP (introduced in \cite{subag:zeitouni:2014}). We shall show that if we consider an SScDPPP on $(0, \infty]$ and take logarithm transform of each atom, then the transformed point process is an SDPPP on $\bar{\real} \setminus \{-\infty\}$. Conversely if we consider an SDPPP on $\bar{\real} \setminus \{-\infty\}$  and take exponential transform of its atoms, then the transformed point process is a SScDPPP point process on $\bar{\real} \setminus [-\infty, 0]$. Based on this connection, we shall derive a slightly extended version (relaxing the property ``increasing" of $h$) of Corollary 3 in \cite{subag:zeitouni:2014}, see Theorem \ref{thm:mainthm_sub_zet} below.

Here, we recall the basic notations and definitions from \cite{subag:zeitouni:2014}. Let $\point =\sum \delta_{p_i}$ be a point process, then by $\theta_x\point$ we denote the shifted point process $\sum \delta_{x+p_i}$ for every $x \in \real$. By $\rmfty$, we denote the space $[-\infty,\infty] \setminus \{- \infty\}$.

\begin{defn}[Decorated Poisson point process, \cite{brunet:derrida:2011}, \cite{maillard:2013}] 
A point process $Q$ is called a decorated Poisson point process of Poisson intensity $m$ and decoration $\point$ (denoted by $Q \sim DPPP(m, \point)$)  if $Q \eql \sum_{i=1}^\infty \theta_{\lambda_i} \point_i$ where $\Lambda = \sum_{i=1}^\infty \delta_{\lambda_i}$ is a Poisson random measure with intensity $m$ and $\point$ is some point process on $\rmfty$ and $\point_i$'s are independent copies of the point process $\point$.
\end{defn}
  
\begin{defn}[Randomly shifted decorated Poisson point process, \cite{subag:zeitouni:2014}]
A point process $T$ is called a \emph{randomly shifted decorated Poisson point process} of Poisson intensity $m$ and decoration $\point$ and shift $U$ (denoted by $T \sim SDPPP(m,\point,U)$) if for $ Q \sim DPPP(m, \point)$ and some independent random variable $U$, $T \eql \theta_U Q$.
\end{defn}

\cite{subag:zeitouni:2014} introduced shift-Laplace functional as
\beq \label{eq:chap2_shift_laplace}
L_T(f|y) = \exptn \bigg( \exp \Big\{ - \int \theta_{-y}f dT \Big\} \bigg)
\eeq
where $\theta_{-y} f(x) = f(x-y)$ for every non-negative measurable function $f : \real \to \real$. By $f \approx g $, we mean the two functions $f$ and $g$ are equal up to translation and let $[g]$ denotes the equivalence class of $g$ under the relation. 

\begin{defn}[Uniquely supported]
A shift-Laplace functional is uniquely supported on $[h]$ if $L_T(f |\cdot) \approx h(\cdot)$ for every $f \in \ckrfty $.
\end{defn}

\begin{thm}[Slight imporovement of Corollary 3 in \cite{subag:zeitouni:2014}] \label{thm:mainthm_sub_zet}
Let $T$ be a point process such that $\prob(T(\rmfty) >0) =1$. Then the following are equivalent.
\begin{enumerate}

\item[(a)] $L_T(f|\cdot)$ is uniquely supported on $[h]$ for some function $h : \real \to \real$.

\item[(b)] $L_T(f |\cdot)$ is uniquely supported on $[h]$, where
\beq \label{eq:chap2_shflaplace}
h(y) = \exptn \bigg( \exp \Big\{ - e^{-c(y-U)} \Big\} \bigg)
\eeq
for some random variable $U$ and $c>0$.

\item[(c)] $T \sim SDPPP(e^{-cx} dx, \point, U)$ for some point process $\point$, random variable $U$ and $c>0$ (same as in (b)).
\end{enumerate}
\end{thm}

\section{ Proof of  ``\eqref{sub:zet:prop3} implies \eqref{sub:zet:prop1}"} \label{sec:slf}

The proof is a heavy-tailed analogue of the proof of Theorem 9 (Converse Part) in \cite{subag:zeitouni:2014}.
Let $N \sim SScDPPP(\malpha, \point, W)$ where $W$ is a positive random variable. We shall compute
\aln{
\Psi_N(f \| y) = \exptn \bigg[  \exptn \bigg(  \exp \big\{ - \sum_{i=1}^\infty \mbfs_W \bbs_{\lambda_i} \point_i (\mbfs_{y^\inv} f) \big\} \bigg| W \bigg) \bigg]  \label{eq:sscdppp_prop31_lap_1}
}
for any $f \in \ckr$ where $\{\lambda_i\}_{i=1}^\infty$ be the atoms of the Poisson point process $\Lambda$ with mean measure $\malpha$ on $(0, \infty]$. Fix $\eta> 0$. Define $I(\eta) = \{ i : \lambda_i > \eta\}$ to be the collection  of all indices of the atoms in the interval $(\eta, \infty]$. We introduce another point process
\alns{
N_\eta = \mbfs_W \sum_{i \in I(\eta)} \mbfs_{\lambda_i} \point_i.
}

It is clear that by monotone convergence theorem, the right hand side of \eqref{eq:sscdppp_prop31_lap_1} equals
\aln{
 \lim_{\eta \to 0} \exptn \bigg[ \exptn \bigg( \exp \bigg\{ - \sum_{i \in I(\eta)} \mbfs_{W \lambda_i} \point_i( \mbfs_{y^\inv} f) \bigg\} \bigg| W \bigg) \bigg]. \label{eq:sscdppp_prop31_lap_2}
}
Note that $|I(\eta)|$ (cardinality of the random set $I(\eta)$) is a Poisson random variable with mean $\eta^{-\alpha}$ and it is independent of $W$ and $\{\point_i\}_{i \ge 1}$. The conditional expectation in \eqref{eq:sscdppp_prop31_lap_2} can then be written as 
\aln{
\exptn \bigg[ \exptn \bigg( \exp \bigg\{ - \sum_{i \in I(\eta)} \mbfs_{W \lambda_i} \point_i(\mbfs_{y^\inv} f) \bigg\} \bigg| I(\eta), W \bigg) \bigg| W \bigg]. \label{eq:sscdppp_prop31_lap_3}
}

It easily follows that conditioned on the event $I(\eta) =k$,
\alns{
\sum_{i \in I(\eta)} \delta_{\lambda_i} \eqd \sum_{i=1}^k \delta_{\eta X_i}
}
for every $\eta > 0$ and $\{X_i\}_{i \ge 1}$ be an i.i.d. collection of  $Pareto(\alpha)$ random variables with probability density function $f_X(x) = \alpha x^{-\alpha -1}$ for $x > 1$. Using this fact, the conditional expectation in \eqref{eq:sscdppp_prop31_lap_3} can be written as,
\aln{
& \exptn \bigg[ \prod_{i \in I(\eta)} \exptn \bigg( \exp \bigg\{ - \point_i (\mbfs_{\eta X_i W y^{\inv}} f) \bigg\} \bigg| I(\eta), W \bigg)  \bigg| W\bigg] \nonumber \\
& = \exptn \bigg[ \bigg( \exptn \Big( \Psi_{\point_1} (f\| y \eta^\inv W^\inv X_1^\inv) \Big) \bigg| W \bigg)^{|I(\eta)|} \bigg| W\bigg]. \label{eq:sscdppp_prop31_lap_4}
}
As $I(\eta)$ is a Poisson random variable with mean $\eta^{-\alpha}$, using the expression for probability generating function for the Poisson random variable we obtain the following expression for the right hand side of \eqref{eq:sscdppp_prop31_lap_4}
\aln{
\exp  \bigg \{ \eta^{-\alpha}  \exptn \Big( \psi_{\point_1} (f \| y \eta^\inv X_1^\inv W^\inv) \Big| W \Big) -1 \bigg\}. \label{eq:sscdppp_prop31_lap_5}
}
It can easily be computed using probability density function of of $X_1$, that 
\aln{
 \exptn \bigg( \Psi_{\point_1} (f \| y \eta^\inv X_1^\inv W ^\inv ) \bigg| W \bigg) 
= y^{-\alpha} \eta^\alpha W^\alpha \int_{x > y^\inv \eta W} \Psi (f\| x^\inv ) \malpha(\dtv x). \label{eq:sscdppp_prop31_exptn_pareto}
}
Combining expressions in \eqref{eq:sscdppp_prop31_exptn_pareto}, \eqref{eq:sscdppp_prop31_lap_5} and \eqref{eq:sscdppp_prop31_lap_2} and using monotone convergence theorem, we get 
\alns{
\Psi_N(f \| y) &=  \exptn \bigg( \exp \bigg\{ - y^{-\alpha} W^{-\alpha} \int (1- \Psi_\point(f \| x^\inv)) \malpha (\dtv x) \bigg\} \bigg) \\
& = \exptn \bigg( \fralpha \Big( y W^\inv c_f \Big) \bigg)
}
where 
\alns{
c_f = \int_{x>0} (1- \Psi_\point(f \| x^\inv) ) \malpha(\dtv x).
}

\section{Properties of the scaled Laplace functional} \label{sec:prop_scaledLaplace}

Suppose that $N$ is point process with scaled Laplace functional which is uniquely supported on $g$. In this section our aim is to study the properties of the function $g$. This study is in parallel to \cite{subag:zeitouni:2014}. We shall show that $g$ is continuous and using continuity of $g$, we derive that $g$ is monotone. Then, we shall show that $g$ is a distribution function which implies that $g$ is increasing. Then using the fact that $N$ is locally finite, we shall show that $(1-g)$ is regularly varying at $\infty$.

Let $N$ be a point process satisfying the assumptions stated in Theorem~\ref{thm:sub_zet}. The first step will be to show that $g$ is a continuous and increasing function. To be more specific, we shall show that $g$ is a distribution function. Then we shall determine the form of $g$. 
 
Continuity of $g$ follows from dominated convergence theorem. Note that
\alns{
 g(y c_f) = \prob(N(\bar{\real}_0) =0) + \prob(N(\bar{\real}_0) >0) \exptn \bigg( \exp \{ - N(\mbfs_{y^\inv} f) \} | N(\bar{\real}_0) >0 \bigg). 
} 
So in is enough to consider the case $\prob(N(\bar{\real}_0) =0) =0$. Following the same same arguments in Lemma 12 of \cite{subag:zeitouni:2014}, it is easy to see that
\aln{
0 = \inf_{y \in \real_+} g(y) < g(x) < \sup_{ y \in \real_+} g(y) =1
}  
and $g$ does not attain its lower bound.

In the next step, we shall show that $g$ is monotone. Let the super-level set of $g$ is denoted by 
\alns{
\suplevel_x = \{ y \in \real_+: g(x) > y\}.
}
To show that $g$ is monotone, it is enough to show that either $\suplevel_x$ has unbounded component or $\suplevel_x$ has component with left end point $0$. Suppose that $\suplevel_x$ has a bounded component which is denoted $(y_x, y'_x)$ such that $y_x > 0$. Following the argument in \cite{subag:zeitouni:2014}, we can show that $\suplevel_x$ has uncountable components. This is contradiction to the fact that $\suplevel_x$ can have at most countable components as $\suplevel_x$ is an open set ($g$ is continuous). Note that we can obtain a relation between components of $\suplevel_x$ and components of super-level set of $\Psi_N(af \| y)$ using scaling relation instead of shift for all $a>0$ and we use the ratio of end points to establish disjointness of the intervals instead of differences.

We shall again use the method of contradiction to show that $g$ does not attain its maximum. Suppose that $g$ attains its maximum and $g$ is increasing. Fix $f \in \ckr$ and define
$$y_0 = \min \{ y \in \real_+: \Psi_N(f \| y) =1\}.$$
Note that $\Psi_N(f \| y) < \Psi_N(2^\inv f \| y)$ for all $y \in \real_+$ i.e. $\Psi_n(2^\inv f \| y) = 1$ for all $y \ge y_0$. Using the fact that $f$ and $2^\inv f$ has the same support, we get that 
\alns{
y_0 = \min \{ y \in \real_+: \Psi_N(2^\inv f \| y) = 1\}
}
Note that $\Psi_N(f \| y)$ is uniquely supported on $g$ for every $f \in \ckr$. This implies that $c_{f} = c_{2^\inv f}$. This means $\Psi_N(f \| y) = \Psi_N(2^\inv \| y)$ for all $y \in \real_+$. This contradicts the fact that $\prob(N(\zrer) =0) =0$.

The following two properties of $g$ can easily be derived in parallel to Corollary 13 and 14 in \cite{subag:zeitouni:2014}. For future reference, we are stating the properties as facts.

\begin{fact} \label{fact:sscdppp_prop12_cont_g}
Suppose that the conditions in Theorem~\ref{thm:sub_zet} holds and $\prob ( N (\bar{\real}_0)=0) =0$. Let $\{c_n\}_{n \ge 1} $ be a sequence of positive real numbers, such that $g(c_n \cdot) \to h(\cdot)$ pointwise, then either $c_n \to 0$ or $\infty$ and $h$ is a constant function with value in $\{0,1\}$, or $c_n \to c $ and $h(y) = g(cy)$.
\end{fact}

\begin{fact} \label{fact:sscdppp_prop12_mable_ext}
Let $f \ge 0$ be a measurable function on $\bar{\real}_0$ and there exists a sequence of functions $\{f_n\}_{n \ge 1}$ in $\ckr$ such that $f_n$ converges to $f$ monotonically and pointwise. Under the assumptions of Fact \ref{fact:sscdppp_prop12_cont_g}, if $\Psi_N(f \| y) \in (0,1)$ for some $y \in \real_+$ then $\Psi_N(f \| \cdot) $ is uniquely supported on the class $[g]_{sc}$.  
\end{fact}

 Fact \ref{fact:sscdppp_prop12_cont_g} easily follows from the continuity of $g$ and  Fact \ref{fact:sscdppp_prop12_mable_ext} follows from a combination of monotone convergence theorem and Fact \ref{fact:sscdppp_prop12_cont_g}.

In this step, we shall show that $g$ is a distribution and hence an increasing function. Define 
\aln{
\maxmod(N) = \inf \{y \in \real_+: N(\{x: |x| > y\}) = 0\} 
}
i.e. maxima of the absolute values of the points associated to the point process $N$. Define $A = \{x \in \bar{\real} : |x|>1\} = [-\infty, 1) \cup (1, \infty]$ and
\alns{
\infty \bbo_A(x) = \begin{cases} \infty & \mbox{ if } x \in A \\ 0 & \mbox{ if } x \notin A. \end{cases}
}
Also note that $\Psi_N(\infty \bbo_A \| y)= \prob(\maxmod(N) \le y)$. It is clear that $\infty \bbo_A \notin \ckr$. Consider the sequence of functions $(f_n : n \ge 1)$ defined below
\begin{align*} 
 f_n(x) = \begin{cases} n  & \mbox{ if } |x| \ge 1+\frac{1}{n}, \\
 n^2 (x-1) & \mbox{ if } 1 \le x < 1+\frac{1}{n},  \\
 n^2 (x+1)  & \mbox{ if } -1-\frac{1}{n} \le x \le -1, \\
 0 & \mbox{ otherwise.}\end{cases}
\end{align*} 
It is easy to see that $f_n \in \ckr$ for all $n \ge 1$ and $f_n \in (0,1)$. Note that $f_n$ converges monotonically and pointwise to $\infty \bbo_A$. So using Fact~\ref{fact:sscdppp_prop12_mable_ext}, we get that there exists a positive real number $\cmax$ which satisfies
\aln{
\prob \Big( \maxmod(N) \le y \Big) = \Psi_N(\infty \bbo_A \| y) = g(\cmax y) \label{eq:sscdppp_prop13_dist_maxmod_gcmax}
}
for all $y \in \real_+$ (in parallel to equation (6.1) in \cite{subag:zeitouni:2014}). Hence $g$ is a distribution function and as a consequence, we show that $g$ is increasing.

In the final step, we study the tail behavior of the distribution function $g$ and show that $1-g$ is regularly varying at $\infty$. Let $\nu$ denotes the mean measure of the point process $N$, that is, $\nu(B) = \exptn (N(B))$ for every Borel subset of $\real$. We shall first show that $\nu((b, \infty)) < \infty$ for all $b>0$. It is clear from  \eqref{eq:chap2_ass_sub_zet} in  Theorem~\ref{thm:sub_zet} that $\nu((ay, \infty])$ is finite for $y \ge 1$. Consider a collection of positive real numbers $(c_a(t): t >0)$ such that $g(c_a(t) y) = \Psi_N(t \bbo_{(a, \infty)} \| y)$. Note that $t \bbo_{(a, \infty)} \notin \ckr$. However, we can construct a sequence of functions $f_n^{(t)} \in \ckr$ such that $f_n^{(t)}$ converges to monotonically and pointwise to the function $t \bbo_{(a, \infty)}$ for every $t>0$. Using Fact~\ref{fact:sscdppp_prop12_mable_ext}, we can show that there exists a constant $c_a(t)$ such that $g(c_a(t) y) = \Psi_N(t \bbo_{(a, \infty)} \| y)$ holds for every $t >0$. In parallel to the steps in Proposition~16 in \cite{subag:zeitouni:2014}, it is easy to see that 
\aln{
\nu((ay, \infty)) = \lim_{t \downarrow 0} \frac{1}{t} \Big(1 - g(y c_a(t)) \Big) < \infty \label{eq:sscdppp_prop12_nu_ass}
}
for all $y \ge 1$. So for all pairs $y_1, y_2 \ge 1$, we have 
\aln{
\frac{\nu( (ay_1, \infty))}{\nu((ay_2, \infty))} = \lim_{t \downarrow 0} \frac{1- g \Big( c_a(t) y_1\Big)}{1- g \Big( c_a(t) y_2 \Big)}. \label{eq:sscdppp_prop12_ratio_nu}
}
It is clear that the right hand side of \eqref{eq:sscdppp_prop12_ratio_nu} is finite as the ratio in left hand side is finite. Using the fact that $g$ is increasing, we get $c_a(t) \uparrow \infty$ as $t \downarrow 0$. So, the right hand side of \eqref{eq:sscdppp_prop12_ratio_nu} becomes
\alns{
\lim_{x \to \infty} \frac{1- g(y_1 y_2^\inv x)}{1- g(x)}
}
and  depends only on the ratio $y_1 y_2^\inv$. This fact implies that 
\aln{
\lim_{x \to \infty} \frac{1-g(y x)}{1- g(x)} \label{eq:sscdppp_prop12_ratio_g}
}
exists and finite for every $y \in \real_+$. Now it is important to note that
\aln{
\nu((ay, \infty)) = \lim_{x \to \infty} \frac{1- g ( y x)}{1- g(x)} \lim_{t\to 0} \frac{1- g \Big( c_a(t) \Big)}{t}.
 \label{eq:sscdppp_prop12_ratio_final}
}
It is clear that the first term in \eqref{eq:sscdppp_prop12_ratio_final} is finite as the ratio \eqref{eq:sscdppp_prop12_ratio_g} is finite for every $y \in \real_+$ and the second term is finite from \eqref{eq:sscdppp_prop12_nu_ass}. Hence we have established the fact that $\nu((b, \infty)) < \infty$  for every $b>0$. Similarly, one can show that $\nu((-\infty, -b)) < \infty$. Finally, these results imply that $\nu(B)< \infty$ for all Borel subsets of $\real$ which are bounded away from $0$.

Let $\theta(y)$ denotes the expression in \eqref{eq:sscdppp_prop12_ratio_g}. Then it is easy to verify that $\theta(\cdot)$ satisfies famous Hamel equation $\theta(yz) = \theta(y) \theta(z)$. So we can write $\theta(y) = y^{-\alpha}$ for all $y \ge 1$ and some $\alpha \in \real$. Using the fact that $g$ is increasing, we get that $\theta(y) \le 1$ for $y \ge 1$. So we obtain $\alpha \ge 0$. Proof of the fact that $\alpha \neq 0$ is very similar to that in Proposition 16 in \cite{subag:zeitouni:2014}. One needs to replace the interval $(0,1)$ by $A$ and shift by scale to prove it. So, we get that there exists  $\alpha >0$, such that 
\alns{
\lim_{x \to \infty} \frac{1- g(xy)}{1-g(x)} = y^{-\alpha} 
} 
for all $y > 0$ (in parallel to equation (6.2) in \cite{subag:zeitouni:2014}). As $g$ is the distribution function of $\maxmod(N)$, it is clear that $\maxmod(N)$ has regularly varying tail at $\infty$.

\section{Distribution of the scale-decoration} \label{sec:scale_decoration}

In this section, we shall construct the scale-decoration $SD(N)$ of $N$ from its scaled Laplace functional. The crucial step will be to compute the distribution of the weak limit $\hat{N}$ of $\mbfs_{y^{\inv}} N$ conditioned on the event $\{\maxmod(N) > y\}$. Then we shall show that the decoration $SD(N)$ corresponding to the point process $N$, has the same distribution as $\mbfs_{(\maxmod(\hat{N}))^\inv}  \hat{N}$.

It is not easy to derive the distribution of $\hat{N}$ directly. We shall construct a collection of point processes $\{N^{(y)} : y \ge 1\}$ from $N$. It will be shown that $\{N^{(y)} : y \ge 1\}$ is a tight family of point processes and we shall obtain the weak limit $N_*$. Then we shall construct another collection of point processes $\{N_{(y)} : y \ge 1\}$. We shall show that the family of the point processes is tight and weak limit $\widetilde{N}$ has the same distribution as $\hat{N}$. From the weak limit $\widetilde{N}$, we shall derive the distribution of $SD(N)$. This construction of scale-decoration is motivated from \cite{subag:zeitouni:2014}. The following property of  $g$ is in parallel to equation (6.3) in \cite{subag:zeitouni:2014}
\alns{
\lim_{m \to 0} c_{mf} = \infty \mbox{ for all } f \in \ckr. 
}

Let $\Novery$ denotes the point process such that for every Borel subset $ B \subset \zrer$, we have
\alns{
\Novery(B) = N \Big( yB \cap \{ (-\infty, -y) \cup ( y, \infty)\} \Big)
} 
conditioned on the event that $\{\maxmod(N) > y\}$. Following the same steps of Lemma 22 in \cite{subag:zeitouni:2014}, it can be shown that $\Novery$ is tight family of point processes. As a consequence, existence of  the limit of the scaled Laplace functionals of $\Novery$ will imply existence of the weak limit $N^*$ of the collection of the point processes $\{\Novery : y \ge 1\}$  as $y \to \infty$. The limit of the scaled Laplace functionals will correspond to the point process $N^*$. Now, we shall derive the scaled Laplace functional of $N_*$ for a function $f \in \ckr$ with support contained in $A$. Exactly the same steps in display (6.4) of \cite{subag:zeitouni:2014}, lead to 
\aln{
 \Psi_{N^*}(f \| x) & = \lim_{y \to \infty} \Psi_{\Novery}(f \|x) \nonumber \\
& = 1 - \lim_{y \to \infty} \frac{1- g(yx c_f \cmax^\inv)}{ 1- g(y)} \nonumber \\
&= 1- x^{-\alpha} ( c_f \cmax^\inv)^{-\alpha}. \label{eq:sscdppp_prop13_lap_nstar}
}
If we consider $f = \infty \bbo_A$ in the  left hand side of \eqref{eq:sscdppp_prop13_lap_nstar}, then in the right hand side, we get $c_f = \cmax$. This implies that $\maxmod{N^*}$ is a $Pareto(\alpha)$ random variable. 

For every $y>0$, we define a point process $\Nundery$ which has same distribution as that of $\mbfs_{y^\inv} N^*$ conditioned on the event $\{\maxmod(N^*) > y\}$. For $y \ge 1$ and any $f \in \ckr$ with support contained in $A$, we have
\aln{
\Psi_{\Nundery}( f) = 1- (c_f \cmax^\inv)^{-\alpha}. \label{eq:sscdppp_prop13_lap_nundery}
}
which is independent of $y$. Define $\Nundery|_A (B) = \Nundery(A \cap B)$. From  \eqref{eq:sscdppp_prop13_lap_nundery}, it is clear that distribution of $\Nundery |_A$ does not depend on $y$. Following the steps of Lemma~24 in \cite{subag:zeitouni:2014}, we obtain
\aln{
\Nundery \eqd N_{(ym)} |_{y^\inv A} \label{eq:sscdppp_prop13_nudery_invar_m}
} 
for all $y \ge 1$ and $m \ge 1$. Our next step is to show that weak limit of $\Nundery$ exists as $y \to \infty$. In order to show existence, we shall show that the family of point processes $\{\Nundery: y \ge 1\}$ is tight and finite dimensional distribution converges. Consider a sequence of real numbers $\{y_n \}$ such that $y_n \to \infty$ as $n \to \infty$. Our aim is to show that for any $f \in \ckr$,
\alns{
\lim_{t \to \infty} \limsup_{n \to \infty} \prob \Big( N_{y_n} (f) >t \Big) = 0.
}
Fix $f \in \ckr$ and choose a large enough $x>1$, such that the support of $f$ is contained in $x^\inv A$ (choice of $x$ depends on $f$). As $y_n \to \infty$, we can find a large enough $n_0$, such that $y_n > x$ for all $n \ge n_0$. For all $n \ge n_0$, we get $y_n = x z_n$ such that $z_n >1$ and hence  we have 
\aln{
N_{(x z_n)}|_{x^\inv A} \eqd N_{(x)}. \label{eq:sscdppp_prop13_tight_nundery}
}
So \eqref{eq:sscdppp_prop13_tight_nundery} follows immediately. The same argument implies convergence of finite dimensional distributions.

Recall that our aim is to find the distribution of $\widehat{N}$ where $\widehat{N}$ is the weak limit of $\mbfs_{y^\inv} N$ conditioned on the event $\{\maxmod(N) > y\}$. We shall show that weak limit $\widetilde{N}$ of  $\Nundery$ has the same distribution as that of $\hat{N}$. Consider a collection of sets $\{A_i\}_{i=1}^l$ which are bounded away from $0$. Then, we get the following expression for distribution of $\widetilde{N}$
\aln{
& \prob \Big( \widetilde{N} ( A_i) \ge k_i , 1 \le i \le l \Big) \nonumber \\
& = \lim_{y \to \infty} \frac{ \prob \Big( N^* (yA_i) \ge k_i , 1 \le i \le l, N^*(yA) >0 \Big) }{\prob(N^*(yA) >0)} \nonumber \\
&= \lim_{y \to \infty}   \lim_{t \to \infty} \frac{\prob(N(tA) >0)}{\prob(N(tyA) >0)}  \frac{\prob \Big( N(ty A_i) \ge k_i, i \le i \le l, N(tyA) >0 \Big)}{\prob( N(tA) >0)}. \label{eq:sscdppp_prop13_eqd_ntilde_nhat}
}
Note that after cancellation, each of the terms in the right hand side of \eqref{eq:sscdppp_prop13_eqd_ntilde_nhat} does not involve $t$ and $y$ separately, but involves the product $ty$. So the separate limits $\lim_{t \to \infty} \lim_{y \to \infty} $ can be replaced by $\lim_{ty \to \infty}$ and we get expression for the distribution of $\widehat{N}$.

Here we shall study some properties of $\widetilde{N}$. As $\Nundery$ converges weakly to $\widetilde{N}$, then $N_{(ym)} |_{y^\inv A}$ converges weakly to $\widetilde{N}|_{y^\inv A}$ as $m \to \infty$. Now use \eqref{eq:sscdppp_prop13_nudery_invar_m} to observe that $\widetilde{N}|_{y^\inv A} \eqd \Nundery$. In parallel to Corollary 26 in \cite{subag:zeitouni:2014}, it is easy to see that for any $f \in \ckr$,
\aln{
\exptn \bigg[ \exp \bigg\{ - \int f \dtv \mbfs_{y^\inv} \widetilde{N} \bigg\} \bigg| \maxmod(\widetilde{N}) > y \bigg] \label{eq:sscdppp_prop13_condinvar_tilden}
}
is independent of $y \ge 1$. It is important to observe that $\maxmod(\widetilde{N}) >1 $. Using the fact that $\widetilde{N} |_A \eqd N^*$, we get that $\maxmod(\widetilde{N}) \eqd \maxmod(N^*)$ and so $\maxmod(\widetilde{N})$ follows $Pareto(\alpha)$ distribution.

Define the scale-decoration point process,
\alns{
SD(N) = \mbfs_{\maxmod(\widetilde{N})^\inv} \widetilde{N}.
}
Our first step will be to observe that $\maxmod(\widetilde{N})$ and $SD(N)$ are independent. In parallel to the steps of Lemma~27 in \cite{subag:zeitouni:2014}, we can write down the conditional Laplace functional of $SD(N)$ given $\{\maxmod(\widetilde{N}) > y\}$ as 
\alns{
\exptn \bigg[ \exp \bigg\{ - \int f \dtv \mbfs_{(\maxmod(\mbfs_{y^\inv} \widetilde{N}))^\inv} \mbfs_{y^\inv} \widetilde{N} \bigg\}  \bigg| \maxmod(N) > y\bigg]
}
which does not involve $y$ using \eqref{eq:sscdppp_prop13_condinvar_tilden}. Hence $SD(N)$ and $\maxmod(\widetilde{N})$ are independent. So we can write $\widetilde{N} = \mbfs_X SD(N)$ where $X$ is a $Pareto(\alpha)$ random variable independent of $SD(N)$.

Let $L(N) \sim ScDPPP(\malpha(\dtv x), \mbfs_{\cmax^\inv} SD(N))$. Then we observe that $L(N)$ has the scaled Laplace functional 
\aln{
\Psi_{L(N)} (f \| y) = \exp \{ - y^{-\alpha} c_f^{-\alpha} \} \label{eq:sscdppp_prop13_scaledLaplace_ln}
}
which is uniquely supported on $\fralpha$. The derivation is similar to the computation done in Section~\ref{sec:slf}.

\section{ Rest of The Proofs} \label{sec:rest_proofs}

\subsection{Proof of  Theorem~\ref{thm:sub_zet}}

Suppose that \eqref{sub:zet:prop2} holds for some positive random variable $W$ and positive scalar $\alpha$. Consider the point process $L'(N) \sim SScDPPP(\malpha(\dtv x), \mbfs_{\cmax^\inv} SD(N), W)$. then in the light of scaled Laplace functional computed in \label{eq:sscdppp_prop13_scaledLaplace_ln}, it follows that the scaled Laplace functional of $L'(N)$ is same as that given in \eqref{sub:zet:prop2}.

Under \eqref{sub:zet:prop1}, we have shown that $g$ is a distribution function. Now we consider a random variable $W_g$ such that $W_g$ follows the distribution $g$ and independent of $L(N) \sim ScDPPP(\malpha(\dtv x),$ $ \mbfs_{\cmax^\inv} SD(N))$. Consider also a F\'{r}echet-$\alpha$ random variable $W_F$ such that $W_F$ is independent of $N$. It is easy to see that
\alns{
\Psi_{\mbfs_{W_g} L(N)} (f) = \Psi_{\mbfs_{W_F} N}(f)
}
using multiplicative convolution. Here it is clear that $\mbfs_{W_F} N \eqd \mbfs_{W_g} L(N)$. Now we shall use the transfer principle to establish that there exists some random variable $\widehat{W_g}$ and a point process $\widehat{L(N)}$ which are independent of each other, such that 
\alns{
\mbfs_{W_F} N = \mbfs_{\widehat{W_g}} \widehat{L(N)}
}
almost surely. Hence we get that $N= \mbfs_{W_F^\inv \widehat{W_g}} \widehat{L(N)}$ establishing \eqref{sub:zet:prop3}. This completes the proof.

\subsection{ Proof of Corollary \ref{propn:strict_stable}}

It is easy to verify that \eqref{item:DM1} implies \eqref{item:DM3} and \eqref{item:DM2} by computing scaled Laplace functional. Suppose that \eqref{item:DM3} holds, then we can construct $SD(N)$ and obtain $\cmax$ so that $L(N) \sim ScDPPP(\malpha(\dtv x), \mbfs_{\cmax^\inv } SD(N))$ has the same scaled Laplace functional as $N$. So \eqref{item:DM1} follows. To show the equivalence, we only have to show that \ref{item:DM3} implies \eqref{item:DM2}.

Using the fact that $N$ is a Radon measure and $\prob(N(\zrer) >0)$, it is clear that $\Psi_N(f \| y) \in (0,1)$ for every $f \in \ckr$ and $y \in \real_+$. Define $\phi(a) = \log \Psi_N (f \| a^\inv)$. Then using the standard way of approximating real numbers by rational ones and then continuity of $\phi(\cdot)$, it is easy to see that 
\alns{
\phi(x) = x^{-\alpha} \phi(1).
} 
Hence it is clear that $\Psi_N(f \| y)$ is uniquely supported on $[\fralpha]_{sc}$.

\subsection{Proof of Theorem \ref{thm:mainthm_sub_zet}}

We shall prove this theorem as a corollary to  Theorem~\ref{thm:sub_zet}. Let $\nmfty$  denotes the space of all counting measures on $\rmfty = \bar{\real} \setminus \{ \infty \}$. We denote the space of all counting measures on $\real_+$ by $\npreal$. Let $\Exp : \nmfty \to \npreal$ be a bijection such that $\Exp(\sum \delta_{a_i}) = \sum \delta_{\exp(a_i)}$ where $a_i \in \rmfty$ for all $i=1,2,\ldots$ and similarly $\Log : \npreal \to \nmfty$ such that $\Log(\sum \delta_{a_i}) = \sum \delta_{\log a_i}$ where $a_i \in \real_+$ for all $i=1,2, \ldots$. It is very easy to see that $\Exp^\inv = \Log$. We also define an operator $\scrt : \ckrfty \to \ckrp$ such that $\scrt(f)(x) = f(\log x)$ for all $x \in \real_+$ and for every $f \in \ckrfty$. It is clear that $\scrt^\inv : \ckrp \to \ckrfty$ such that $\scrt^\inv (u)(x) = u (e^x)$ for all $x \in \rmfty$ and for every $u \in \ckrp $. Now we would like to state some easy consequences of change of variable formula. Suppose that $T \in \nmfty$ and $u \in \ckrp$ and $t>0$, then
\beq
\int \mbfs_t u \dtv \Exp(T) = \int \theta_{\log t} \scrt^\inv(u) \dtv T \label{eq:cov1}
\eeq
similarly, we get
\beq \label{eq:cov2}
\int \mbfs_t u \dtv N = \int \theta_{\log t} \scrt^\inv(u) \dtv \Log N
\eeq
where $N \in \npreal$ and $u \in \ckrp$. Suppose that $T \in \nmfty$ and shift-Laplace functional of $T$ is shift-uniquely supported. Then our first step is to show that scaled Laplace functional of $\Exp T$ is scale-uniquely supported. Consider $h \in \ckrp$ and $t >0$ to see that
\begin{align*}
\Psi_{\Exp(T)} (\mbfs_{t^\inv} u) &= \exptn \bigg( \exp \Big\{ - \int \mbfs_{t^\inv} u \dtv \Exp(T) \Big\} \bigg)   \\
& = h \Big(\log t - \tau_{\scrt^\inv(u)} \Big) \\
& = g \Big(t e^{-\tau_{\scrt^\inv(u)}}\Big)
\end{align*}
where $g(x) = h( \log x)$. Hence we are done with the first step.

Our next step will be to show that if $N \sim SScDPPP(\malpha, \point, W)$ then $\Log(N)\sim  SDPPP(e^{-\alpha x}\dtv x, $ $ \Log(\point), \log W)$. Suppose that our claim is true. Then using the first step and Theorem~\ref{thm:sub_zet}, we obtain a point process $N$ admitting SScDPPP representation. Now we use the claim to observe that $\Log N$ is an SSDPPP and shift-uniquely supported.

So we are only remained with proof of the claim. We shall prove the claim by computing Laplace functional of $\Log N$. Fix $f \in \ckrfty$. We get
\begin{align*}
\Psi_{\Log(N) }(f) & = \exptn \bigg( \exp \Big\{ - \int f \dtv \Log(N) \Big\} \bigg) \\
& = \exptn \bigg( \exp \Big\{ - \int \scrt(f) \dtv N \Big\} \bigg) \\
& = \exptn \bigg( \exp \Big\{-  W^\alpha c_{\scrt(f)}^{-\alpha}   \Big\} \bigg)\\
&= \exptn \bigg( \exp \Big\{ - e^{\alpha(\log W- \log c_{\scrt(f)})} \Big\} \bigg)
\end{align*}
Now to get the complete description of the decoration we need to compute
\begin{align*}
c_{\scrt(f)} &= \int_0^\infty \Big(  1- \Psi_{\point}(\mbfs_x \scrt(f)) \Big) \alpha x^{-\alpha -1} \dtv x \\
&= \int_{-\infty}^\infty e^{-\alpha x} \bigg( 1- L_{\Log(\point)}(f|-x) \bigg) \dtv x
\end{align*}
Finally we get,
\begin{align*}
\Psi_{\Log(N) }(f) =  \exptn \bigg( \exp \Big\{ - e^{\alpha \Big[\log W- \log \Big( \int_{-\infty}^\infty e^{-\alpha x} \big( 1- L_{\Log(\point)}(f|-x) \big) \dtv x\Big)   \Big]} \Big\} \bigg)
\end{align*}
which implies that $ \Log N \sim SSDPPP(e^{-\alpha x} \dtv x, \Log(\point), \log W)$.

\section*{Acknowledgement}
The research was supported by NWO VICI grant. The author is thankful to Rajat Subhra Hazra and Parthanil Roy for numerous discussions and helpful suggestions which improved   this article.  

\small
\bibliographystyle{abbrvnat}

\vspace{.3cm}

\noindent \footnotesize{Centrum Wiskunde \& Informatica \\P.O. Box 94079\\
1090 GB Amsterdam\\
NETHERLANDS.\\
Email address: ayan.bhattacharya@cwi.nl}

\end{document}